\documentclass[11pt, a4paper]{amsart}
      \usepackage{amssymb}

      \theoremstyle{plain}
      \newtheorem{theorem}{Theorem}[section]
      \newtheorem{lemma}[theorem]{Lemma}
      \newtheorem{corollary}[theorem]{Corollary}
      \newtheorem{proposition}[theorem]{Proposition}

      \theoremstyle{definition}
      \newtheorem{definition}[theorem]{Definition}

      \theoremstyle{remark}

      \makeatletter
      \def\@setcopyright{}
      \def\serieslogo@{}
      \makeatother

\begin{document}

% First we specify the top matter (author, title, etc).
%
% Note: All of the top matter items are optional and can be omitted.
% But you will probably want to specify at least the author and title
% and perhaps an abstract.

   % author information

   % first author

   \author{Sungwoon Kim}
   \address{School of Mathematics,
   KIAS, Heogiro 85, Dongdaemun-gu,
   Seoul, 130-722, Republic of Korea}
   \email{sungwoon@kias.re.kr}

   % second author

   \author{Inkang Kim}
   \address{School of Mathematics,
   KIAS, Heogiro 85, Dongdaemun-gu
   Seoul, 130-722, Republic of Korea}
   \email{inkang@kias.re.kr}

   % the address where the research was carried out
   %\address{}

   % current address, usually not needed because it is the same as the
   % regular address
   %\curraddr{}

   %\email{}

   % title

   \title[Simplicial volume of $\mathbb{Q}$-rank one manifolds]{Simplicial volume of
   $\mathbb{Q}$-rank one locally symmetric manifolds covered by the product of
   $\mathbb{R}$-rank one symmetric spaces}

   % Note that the short title for running heads goes in square
   % brackets.  This is optional.  The long title goes in curly
   % braces.  In the long title, line breaks are indicated by \\.

   % abstract (optional)
\begin{abstract}
    In this paper, we show that the simplicial volume of $\mathbb{Q}$-rank one
    locally symmetric spaces covered by the product of $\mathbb{R}$-rank one symmetric
    spaces is strictly positive.
\end{abstract}

 \footnotetext[2]{The second
 author gratefully acknowledges the partial support of KRF grant
(0409-20060066).}

   % AMS subject classifications (used in AMS journals)

   % AMS keywords (used in AMS journals)

   % acknowledge support, etc
   \thanks{}
   \thanks{}

   % dedication
   \dedicatory{}

   % today's date, or fill in whatever date you prefer
   \date{}

% This ends the top matter information.
% We can now tell LaTeX to display the top matter.

   \maketitle

% Having displayed the top matter, we now proceed to the body of the
% article.

% The body of the article is divided into sections.
% Each section begins with a \section command.

    \section{Introduction}

The simplicial volume of a connected, oriented manifold $M$ was introduced by Gromov \cite{Gro82}.
This is a topological invariant in $\mathbb{R}_{\geq 0}$ and measures how
efficiently the fundamental class of $M$ can be represented by simplices.
Gromov conjectured that non-positively curved closed manifolds with negative Ricci curvature
have positive simplicial volume.

First, the positivity of the simplicial volume was verified for closed negatively curved manifolds by Thurston \cite{Th78}, Gromov \cite{Gro82} and Inoue-Yano \cite{IY82}.
It was verified for closed locally symmetric spaces covered by $SL_3(\mathbb{R})/SO_3(\mathbb{R})$ by
Savage \cite{Sa82} and Bucher-Karlsson \cite{BK07}.
Then, Lafont and Schmidt \cite{LS06} showed that the simplicial volume of all closed locally symmetric spaces of non-compact type is positive, which gave a positive answer to the conjecture raised by Gromov.

For open manifolds, the simplicial volume is somewhat mysterious.
Thurston \cite{Th78} verified that the simplicial volume of complete
Riemannian manifolds with pinched negative sectional curvature and
finite volume is strictly positive. In contrast, Gromov
\cite{Gro82}, L\"{o}h and Sauer \cite{LS09-2} proved that the
simplicial volume of open manifolds, which are the Cartesian product
of three open manifolds and locally symmetric spaces of
$\mathbb{Q}$-rank  at least $3$, vanishes. L\"{o}h and Sauer
\cite{LS09-1} showed that Hilbert modular varieties have positive
simplicial volume, which was the first class of examples of open
locally symmetric spaces of $\mathbb{R}$-rank at least $2$ for which
the positivity of simplicial or minimal volume is known. Hilbert
modular varieties are special cases of $\mathbb{Q}$-rank one locally
symmetric spaces covered by the product of hyperbolic planes. The
aim of this paper is to show the positivity of the simplicial volume
of $\mathbb{Q}$-rank one locally symmetric spaces covered by the
product of $\mathbb{R}$-rank one symmetric spaces.

\begin{theorem} \label{thm1}
Let $M$ be a $\mathbb{Q}$-rank one locally symmetric
space covered by the product of $\mathbb{R}$-rank one symmetric
spaces. Then, the simplicial volume of $M$ is positive.
\end{theorem}

For general $\mathbb{Q}$-rank one locally symmetric spaces,
see the forthcoming paper \cite{BKK}.
Gromov \cite{Gro82} proved a lower bound for the minimal volume of
$n$-dimensional Riemannian manifolds in terms of the simplicial
volume: $$||M||\leq (n-1)^n n!\cdot \text{minvol}(M).$$ The theorem
implies the positivity of minimal volume of $\mathbb{Q}$-rank one
locally symmetric spaces covered by the product of $\mathbb{R}$-rank
one symmetric spaces as follows. See Connell and Farb \cite{CF1} for different approach using Lipschitz
class of locally symmetric metrics.

\begin{corollary}\label{minvolume}
 The minimal volume of
$\mathbb{Q}$-rank one finite volume locally symmetric spaces covered
by the product of $\mathbb{R}$-rank one symmetric spaces is
positive.
\end{corollary}

So far, the degree theorem for the open locally symmetric spaces of non-compact type with
finite volume holds with the Lipschitz condition on the map.
From Theorem \ref{thm1}, one can obtain the degree theorem without
any Lipschitz condition on the map.

\begin{theorem}
Let $N$ be a Riemannian $n$-dimensional manifold of finite volume with Ricci curvature
bounded below by $-(n-1)$ and $M$ be a $\mathbb{Q}$-rank one locally
symmetric space covered by the product of $\mathbb{R}$-rank one symmetric spaces.
For any proper map $f:N\rightarrow
M$ we have
$$deg(f)\leq C_n \frac{vol(N)}{vol(M)},$$ where $C_n$ depends only on $n$.
\end{theorem}

The degree theorem for general locally symmetric space of finite
volume is proved by Connell and Farb \cite{CF1, CF2} with Lipschitz condition on $f$. The essential part of our approach is to show that the
geodesic straightening map is well-defined on the locally finite
chain complex of $\mathbb{Q}$-rank one locally symmetric spaces. In
fact, Thurston \cite{Th78} introduced the geodesic straightening map
on the singular chain complex of non-positively curved manifolds,
which is homotopic to the identity.

Unfortunately, the geodesic straightening map is generally not
defined on the locally finite chain complex of non-positively curved
manifolds because the geodesic straightening of a locally finite
chain is not necessarily locally finite. However, the situation in
the $\mathbb{Q}$-rank one locally symmetric spaces is different. By
using Leuzinger's explicit geometric description of
$\mathbb{Q}$-rank one locally symmetric spaces \cite{Le95}, one can
see that the geodesic straightening of a locally finite chain  is
locally finite. The presence of the geodesic straightening map on a
locally finite chain complex and the uniform upper bound of the
volume of geodesic simplices in $\mathbb{Q}$-rank one locally
symmetric spaces covered by the product of $\mathbb{R}$-rank one
symmetric spaces give rise to the positivity of the simplicial
volume.

% Sections can be labeled for cross referencing.

\section{Preliminaries}

    In this section, we first collect some definitions and results about the simplicial volume.
    We begin with the definition of the simplicial volume.

    \subsection{Simplicial volume}

        Let $M$ be an $n$-dimensional, connected, oriented manifold. Denote by $C_*(M)$
        the singular chain complex of $M$ with real coefficients. Consider on $C_*(M)$ the $\ell^1$-norm
        with respect to the canonical basis of singular simplices, that is,
        $\| c \| _1=\sum_{i=1}^r | a_i |$
        for $c=\sum_{i=1}^r a_i \sigma_i$ in $C_*(M)$. This norm induces a semi-norm on
        the real coefficient homology $H_*(M)$ of $M$ as follows:
        $$ \| z \|= \inf_c  \| c\|_1,$$
        where $c$ runs over all singular cycles representing $z\in H_*(M).$

        The simplicial volume $\| M \|$ of a closed manifold $M$ is defined as the
        semi-norm of the fundamental class $[M]\in H_n(M)$.

        For closed Riemannian manifolds, Gromov proved the remarkable proportionality principle
        relating the simplicial volume and the volume of the Riemannian manifolds \cite{Gro82}.

        \begin{theorem}[Gromov]
        Let $M$ and $N$ be two closed Riemannian manifolds with isometric universal covers. Then,
        $$ \frac{\| M \|}{vol(M)}= \frac{\| N \|}{vol(N)}.$$
        \end{theorem}

        If $M$ is open,
        one cannot choose the fundamental class of $M$ in the singular homology of $M$
        because the top dimensional homology of $M$ is trivial.
        Hence, the simplicial volume of open manifolds is defined in terms of the locally finite chain complex of the manifolds as follows.

        Let $M$ be a topological space and let $S_k(M)$ be the set of all continuous maps
        from the standard $k$-simplex $\Delta^k$ to $M$. A subset $A$ of $S_k(M)$ is called
        \textit{locally finite} if any compact subset of $M$ intersects the image of only a finite number of elements of $A$. Let us denote by $S^\text{lf}_k(M)$ the set of all locally finite subsets of $S_k(M)$.

        \begin{definition}
        Let $M$ be a topological space and let $k\in\mathbb{N}$. The locally finite chain
        complex of $M$ is the chain complex $C^\text{lf}_*(M)$ consisting of the real vector spaces
        $$ C^{\text{lf}}_k(X)=\bigg\{ \sum_{\sigma \in A} a_\sigma\cdot\sigma ~\Big|~ A\in
        S^{\text{lf}}_k(X) \text{ and }(a_\sigma)_{\sigma \in A} \subset \mathbb{R} \bigg\}$$
        equipped with the boundary operator given by the alternating sums of the
        $(k-1)$-faces. The locally finite homology
        $H^\text{lf}_*(M)$ of $M$ is the homology of the locally finite chain complex $C^\text{lf}_*(M)$.
        \end{definition}

        The $\ell^1$-norm $\| \cdot \|_1$ on the locally finite chain complex of $M$ is defined with respect to the canonical basis of singular simplices, and
        it also gives rise to a semi-norm on the locally finite homology of $M$.
        Any oriented, connected manifold $M$ possesses a fundamental class, which is a distinguished
        generator of the top dimensional locally finite homology $H^{\text{lf}}_n(M;\mathbb{Z})\cong\mathbb{Z}$
        with integral coefficients \cite{LS07}. Now, we are ready to define the simplicial volume
        of open manifolds.

        \begin{definition}
        Let $M$ be a connected $n$-dimensional manifold without boundary.
        Then, the simplical volume of $M$ is defined as
        $$\| M\|=\inf\{\| c \|_1~\big|~c\in C^{\text{lf}}_n(M)\text{ is a fundamental cycle of } M\}.$$
        \end{definition}

        The simplicial volume of open manifolds is zero in a large number of cases, including
        the product of three open manifolds \cite{Gro82},
        locally symmetric manifolds of $\mathbb{Q}$-rank of at least $3$ \cite{LS09-2}.

        Gromov \cite{Gro82} introduced another notion of the  simplicial volume, so called Lipschitz simplicial volume which is a
         geometric variant of the ordinary simplicial volume.
        Let $M$ be an oriented Riemannian manifold.
        For a locally finite chain $c=\sum_{\sigma \in A} a_\sigma \sigma$, define
        Lip$(c)$ as the supremum of the Lipschitz constants of the singular simplices $\sigma$
        with respect to the standard Euclidean metric on the standard simplex.

        \begin{definition}
        Let $M$ be an $n$-dimensional, oriented Riemannian manifold.
        A locally finite chain $c$ is called a Lipschitz fundamental cycle of $M$ when $c$
        represents the fundamental class of $M$ and Lip$(c)<\infty$.
        The Lipschitz simplicial volume $\| M\|_\text{Lip}$ of $M$ is defined as
        $$\| M\|_\text{Lip}=\inf\{\| c \|_1~\big|~c\in C^{\text{lf}}_n(M)\text{ is a
        Lipschitz fundamental cycle of } M \}.$$
        \end{definition}

        From the definition of Lipschitz simplicial volume, we have the obvious inequality $\| M\| \leq \| M\|_\text{Lip}$
        for an oriented, Riemannian manifold $M$. If $M$ is closed, the fundamental
        cycles in the locally finite chain complex of $M$ involve only a finite number of simplices, and hence,
        $\| M\| = \| M\|_\text{Lip}$. L\"{o}h and Sauer \cite{LS09-2} prove the proportionality
        principle for the Lipschitz simplicial volume under the non-positive curvature condition
        in the non-compact case.

        \begin{theorem}[L\"{o}h and Sauer] \label{lpp}
        Let $M$ and $N$ be complete, non-positively curved Riemannian manifolds of
        finite volume. Assume that their universal covers are isometric. Then,
        $$ \frac{\| M \|_\text{Lip}}{vol(M)}=\frac{\| N \|_\text{Lip}}{vol(N)}.$$
        \end{theorem}

By Theorem \ref{lpp}, they can show that the Lipschitz simplicial volume of locally
symmetric spaces of finite volume and non-compact type is positive and, moreover,
they obtain degree theorems for locally symmetric spaces of non-compact type of
finite volume, which is originally due to \cite{CF2}.

%\begin{theorem}[Connell and Farb]\label{degth}
%For every $n\in \mathbb{N}$, there is a constant $C_n > 0$ with the following property:
%Let $M$ be an $n$-dimensional locally symmetric space of non-compact type with finite volume.
%Let $N$ be an $n$-dimensional complete Riemannian manifold of finite volume with
%Ricci($N$)$\geq -(n-1)$ and sec($N$)$\leq 1$, and let $f : N \rightarrow M$
%be a proper Lipschitz map. Then
%$$ \text{deg}(f)\leq C_n \frac{\text{Vol}(N)}{\text{Vol}(M)}.$$
%\end{theorem}

\begin{theorem}[Connell and Farb]\label{degth}
Let $M$ be a locally symmetric $n$-manifold of non-compact type with finite volume.
Assume that $M$ has no local direct factors locally isometric to $\mathbb R$, $\mathbb{H}^2$ or $SL_3(\mathbb{R})/SO_3(\mathbb{R})$.
Then for any complete Riemannian manifold $N$ with finite volume and any proper Lipschitz map $f:N \rightarrow M$,
$$deg(f) \leq C \frac{vol(N)}{vol(M)}$$
where $C>0$ depends only on $n$ and the smallest Ricci curvatures of $N$ and $M$.
\end{theorem}

Note that the results of Connell and Farb, L\"{o}h and Sauer hold with the Lipschitz condition on $f$.
If the positivity of the ordinary simplicial volume of $M$ is verified,
one can obtain the degree theorem without the Lipschitz condition on $f$.

\section{$\mathbb{Q}$-rank one locally symmetric spaces}

    In this section, we recall the definitions of arithmetic lattices, $\mathbb{Q}$-rank, and cusp
    decomposition in $\mathbb{Q}$-rank one locally symmetric spaces.
    The cusp decomposition in quotient manifolds by arithmetic lattices is crucial to show the presence of the geodesic straightening map on the locally finite chain complex.

    Let $X$ be a connected symmetric space of non-compact type. Let $G$ be the
    the identity component of the isometry group of $X$. Then, $G$
    is a connected semi-simple Lie group with trivial center and no compact factor
    \cite{Eb96}. We first recall the definition of arithmetic lattices in \cite{Zi84}.

    \begin{definition}
    Let $G$ be a connected semi-simple Lie group with trivial center and no compact factors.
    Let $\Gamma \subset G$ be a lattice. Then, $\Gamma$ is called \textit{arithmetic} if there exist
    \begin{itemize}
    \item[(i)] a semi-simple algebraic group $\mathbf{G}\subset
        GL(n,\mathbb{C})$ defined over $\mathbb{Q}$ and
    \item[(ii)]  a surjective homomorphism $\rho :\mathbf{G}(\mathbb{R})^0
        \rightarrow G$ with compact kernel
    \end{itemize}
    such that $\rho(\mathbf{G}(\mathbb{Z}) \cap
    \mathbf{G}(\mathbb{R})^0)$ and $\Gamma$ are commensurable.
    \end{definition}

    The $\mathbb{Q}$-rank($\Gamma$) is defined as the dimension of any
    maximal $\mathbb{Q}$-split
    torus of $\mathbf{G}(\mathbb{Q})$ when $\Gamma$ is an arithmetic lattice.
    The structure of the ends of $M=\Gamma \backslash X$ is closely
    related to the $\mathbb{Q}$-rank($\Gamma$). For instance, a locally symmetric space
    $M=\Gamma \backslash X$ is compact if and
    only if the $\mathbb{Q}$-rank($\Gamma$) is zero by the result of Borel and Harish-Chandra.
    To understand the ends of quotient manifolds by arithmetic lattices, we recall the reduction theory due
    to A. Borel and Harish-Chandra \cite{Bo69}.

    \begin{theorem}[Borel, Harish-Chandra] \label{prt}
    Let $\mathbf{G}$ be a semi-simple algebraic group defined over $\mathbb{Q}$ with
    associated Riemannian symmetric space $X$. Let $\mathbf{P}$ be a minimal
    parabolic $\mathbb{Q}$-subgroup of $\mathbf{G}$ and let $\Gamma$ be an arithmetic
    subgroup of $\mathbf{G}(\mathbb{Q})$. Then
    \begin{itemize}
    \item[(i)] the set of double cosets $\mathcal{F}=\Gamma \backslash
        \mathbf{G}(\mathbb{Q})/ \mathbf{P}(\mathbb{Q})$ is finite,
    \item[(ii)]  there exists a generalized Siegel set $\mathcal{S}_{\omega
        ,\tau}$ such that for a (fixed) set $\{ q_i\text{ }|\text{ }1\leq i \leq m
        \}$ of representatives of $\mathcal{F}$ the union
        $\Omega=\bigcup_{i=1}^m q_i \cdot \mathcal{S}_{\omega,\tau}$ is a fundamental set
        for $\Gamma$ in $X$.
    \end{itemize}
    \end{theorem}

    By using the reduction theory, Leuzinger \cite{Le95} gives an explicit differential geometric description of $M$.

    \begin{theorem}[Leuzinger]
    Let $X$ be a Riemannian symmetric space of non-compact type
    and with $\mathbb{R}$-rank $\geq 2$ and let $\Gamma$ be an irreducible, torsion-free,
    non-uniform lattice in the isometry group of $X$. On the locally symmetric space $M=\Gamma
    \backslash X$ there exists a continuous and piecewise real analytic exhaustion function
    $h:M\rightarrow [0,\infty)$ such that, for any $s\geq 0$, the sublevel set $M(s)=\{ h\leq s \}$ is
    a compact submanifold with corners of $M$. Moreover the boundary of
    $M(s)$, which is a level set of $h$, consists of projections of subsets of horospheres in $X$.
    \end{theorem}

    More precisely, there exists a union of a countable number of open horoballs
    $U(s)$ of $X$ such that $X(s)=X-U(s)$ is $\Gamma$-invariant and $M(s)=\Gamma \backslash X(s)$.
    If $s_1 < s_2$, then $U(s_2) \subset U(s_1)$ and the subsets $X(s)=X-U(s)$
    exhaust $X$. These horoballs are in one-to-one correspondence with the vertices of the Tits
    building of $\mathbf{G}$ over $\mathbb{Q}$. The deleted horoballs are disjoint if and only if $\Gamma$ is an arithmetic
    subgroup of a semisimple algebraic group of $\mathbb{Q}$-rank one. In the case of higher
    $\mathbb{Q}$-rank, the horoballs of $U(s)$ intersect and give rise to corners \cite{Le04}.

    Indeed, the theorem of Leuzinger is available for torsion-free arithmetic
    lattices with $\mathbb{Q}$-rank at least one because the proof in the paper \cite{Le95} is only based on the
    reduction theory for arithmetic lattices. He used the arithmeticity theorem of Margulis to show
    that all irreducible, torsion-free, non-uniform lattices with $\mathbb{R}$-rank$\geq 2$ are arithmetic
    lattices.

    \section{Geodesic straightening map on locally finite chain complex}

    The geodesic straightening map has played an important role in proving the positivity of
    the simplicial volume. In this section, we will show the presence of the geodesic
    straightening map on the locally finite chain complex of $\mathbb{Q}$-rank one
    locally symmetric spaces.

    \subsection{Geodesic Straightening}

        The geodesic straightening map on the level of singular chain complexes was
        introduced by Thurston \cite{Th78}. We recall the definition of the
        geodesic straightening map.

        Let $X$ be a simply connected, complete Riemannian manifold with
        non-positive sectional curvature. For $x_0,\ldots,x_k \in X$, the geodesic simplex
        $[x_0,\ldots,x_k]$ is defined inductively as follows.
        First, $[x_0]$ is the point $x_0\in X$, and $[x_0,x_1]$
        is the unique geodesic arc from $x_1$ to $x_0$. In general, $[x_0,\ldots,x_k]$
        is the geodesic cone on $[x_0,\ldots,x_{k-1}]$ with the top point $x_k$.

        \begin{definition}
        Let $M$ be a connected, complete Riemannian manifold with non-positive sectional
        curvature.
        Then, the geodesic straightening map $st_*:C_*(M)\rightarrow C_*(M)$ is defined by
        $$st_k(\sigma)=\pi_M\circ [\tilde{\sigma}(e_0),\ldots,\tilde{\sigma}(e_k)]
        \text{ for a singular }k\text{-simplex } \sigma ,$$
        where $\pi_M :\widetilde{M}\rightarrow M$ is the universal covering,
        $e_0,\ldots,e_k $ are the vertices of the standard $k$-simplex $\Delta^k$, and
        $\widetilde{\sigma}$ is a lift of $\sigma$ to the universal cover $\widetilde{M}$.
        \end{definition}

        The following proposition proved by Thurston makes it possible to obtain
        the simplicial volume of $M$ by considering only the $\ell^1$-norm on the geodesically
        straight chains of $M$.

        \begin{proposition}[Thurston]
        Let $M$ be a connected, complete Riemannian manifold with non-positive
        sectional curvature. Then, the geodesic straightening map is chain homotopic
        to the identity.
        \end{proposition}

\subsection{Straightening locally finite chains}

Let $M$ be a $\mathbb{Q}$-rank one locally symmetric space. First, we fix
some notations. Let $X$ denote the universal cover of $M$ and $\Gamma$ the
fundamental group of $M$. Let $h:M\rightarrow [0,\infty)$ be
the exhaustion function of $M$ under Leuzinger's theorem in \cite{Le95}. For
any $s\geq 0$, $M$ admits the following disjoint decomposition
$$ M=M(s) \cup \coprod_{i=1}^l E_i(s),$$
where the sublevel set $M(s)=\{h\leq x \}$ is a compact submanifold and
$E_i(s)$ is a cusp end of $M-M(s)$ for each $i=1,\ldots,l$.
As we already mentioned, there is a countable number of
pairwise disjointed horoballs $U(s)$ in $X$ such that each $E_i(s)$ is obtained by
the quotient of an open horoball in $U(s)$ by $\Gamma$ and $M(s)=\Gamma
\backslash (X-U(s))$. Thus, $E_i(s)$ is geodesically convex for each $i=1,\ldots,l$.

\begin{lemma} \label{str}
Let $M$ be a $\mathbb{Q}$-rank one locally symmetric space.
Then, the geodesic straightening of a locally finite chain in $C^\text{lf}_*(M)$
is a locally finite chain.
\end{lemma}

\begin{proof}
Let $A \in S^\text{lf}_k(M)$ and $c = \sum_{\sigma \in A} a_\sigma \sigma$
be a locally finite chain in $C^\text{lf}_k(M)$. Let us define $st_k(A) = \{
st_k(\sigma)\text{ }|\text{ }\sigma \in A  \}$. It is sufficient to prove that
$st_k(A)$ is locally finite. Let $K$ be a compact subset of $M$.
Then, one can choose a compact submanifold $M(s)$ of $M$ containing $K$ for
some $s>0$.

By the local finiteness of $A$, a compact submanifold $M(s)$
intersects the image of only a finite number of elements of $A$. Let
$\sigma$ be an element of $A$ whose image does not intersect $M(s)$.
Then, we claim that the image of $\sigma$ has to be contained in
only one cusp end of $M-M(s)$. Suppose the image of $\sigma$
intersects at least two cusp ends of $M-M(s)$, denoted by $E_1(s)$
and $E_2(s)$. Since the image of $\sigma$ is path-connected, there
is a path in the image of $\sigma$ connecting two different points
contained in $E_1(s)$ and $E_2(s)$, respectively. However, any path
connecting such two points must pass through $M(s)$. This means that
the image of $\sigma$ intersects $M(s)$, which contradicts the
assumption that the image of $\sigma$ does not intersect $M(s)$.

Now, let's assume that
the image of $\sigma$ is contained in $E_1(s)$.
Since $E_1(s)$ is geodesically convex and the image of $\sigma$ is contained in $E_1(s)$,
the image of geodesic straightening $st_k(\sigma)$ of $\sigma$ is also totally contained
in $E_1(s)$.
This implies that the image of $st_k(\sigma)$ does not intersect both $M(s)$ and $K$.
Hence, we conclude that $K$ can intersect the image of $st_k(\tau)$ for
only a finite number of elements $\tau$ of $A$ intersecting $M(s)$, and so
$st_k(A)$ is locally finite. This completes the proof of the lemma.
\end{proof}

By Lemma \ref{str}, the geodesic straightening map is well-defined
on the locally finite chain complex of $M$: $$st^\text{lf}_* :
C^{\text{lf}}_*(M)\rightarrow C^{\text{lf}}_*(M).$$ The map $st^\text{lf}_*$
is obviously a chain map because it
is induced from the geodesic straightening map on the singular chain
complex of $M$. Furthermore, we prove that it is chain homotopic to the
identity as follows.

\begin{proposition} \label{ch}
Let $M$ be a $\mathbb{Q}$-rank one locally symmetric space.
Then the geodesic straightening map $st^\text{lf}_*$ is chain homotopic
to the identity.
\end{proposition}

\begin{proof}
First, we recall the construction of the chain homotopy
$H_*:C_*(M)\rightarrow C_{*+1}(M)$ from the geodesic straightening map
to the identity. The chain homotopy $H_k$ is defined by the straight line
homotopy between any $k$-simplex and its geodesically straight
simplex. Moreover, these homotopies, when restricted to lower dimentional
faces, agree with the homotopies canonically defined on those faces. For more
details, let $H_\sigma$ be a canonical straight line homotopy $$ H_\sigma :
\Delta^k \times [0,1] \rightarrow M$$ from $\sigma$ to $st_k(\sigma)$ for
any $k$-simplex $\sigma$. Now $\Delta^k \times [0,1]$ has vertices
$$a_0=(e_0,0),\ldots, a_k=(e_k,0),b_0=(e_0,1),\ldots,b_k=(e_k,1).$$
For each $i=0,\ldots,k,$ let
$$\alpha_i : \Delta^{k+1}\rightarrow \Delta^k \times [0,1]$$ be the affine map
that maps $e_0,\ldots,e_{k+1}$ to $a_0,\ldots,a_i,b_i,\ldots,b_k,$
respectively. Define linear transformation $$H_k :
C_k(M;\mathbb{R})\rightarrow C_{k+1}(M;\mathbb{R})$$ by the formula
$$H_k(\sigma)=\sum^k_{i=0}(-1)^iH_\sigma \circ\alpha_i.$$ This $H_*$ is
the chain homotopy from the geodesic straightening map to the identity on
the singular chain complex of $M$.

Let $c=\sum_{\sigma \in A} a_\sigma \sigma$ be a locally finite chain. We
claim that $H_*(c)$ is a locally finite chain again. Let $K$ be a compact subset
of $M$. Choose a compact submanifold $M(s)$ containing $K$ for some $s>0$.
Suppose that the image of $\sigma \in A$ does not intersect $M(s)$.
By a similar argument in Lemma \ref{str},
the images of both $\sigma$ and $st_k(\sigma)$ are totally contained in
only one cusp end, denoted by $E(s)$.

Since $E(s)$ is geodesically convex, the image of
the straight line homotopy $H_k(\sigma)$ between $\sigma$ and
$st_k(\sigma)$ is totally contained in $E(s)$. This means that $M(s)$ does not
intersect the image of $H_k(\sigma)$ and $K$ does not either. Therefore, $K$
can intersect the image of $H_k(\sigma)$ for only a finite number of elements
$\sigma$ of $A$ intersecting $M(s)$. Since $H_k(\sigma)$ is a finite sum of
simplices, $K$ intersects the image of only a finite number of simplices in $H_k(c)$.
In other words, $H_k(c)$ is a locally finite chain. Finally, we obtain
the following well-defined map:
$$H^\text{lf}_k : C^\text{lf}_k(M)\rightarrow C^\text{lf}_{k+1}(M)$$
satisfying
$$\partial_{k+1}H^\text{lf}_k+H^\text{lf}_{k-1}\partial_k= st^\text{lf}_k-id.$$
Therefore, $H_*^\text{lf}$ is a chain homotopy from $st^\text{lf}_*$ to the identity.
\end{proof}

As can be seen in the proof of Lemma \ref{str} and Proposition
\ref{ch}, $\mathbb{Q}$-rank one condition on $M$ is essential to
obtain the geodesic straightening map on the locally finite chain
complex of $M$. Since the cusp end of higher $\mathbb{Q}$-rank
locally symmetric space has corners, it is not geodesically convex.
Hence, Lemma \ref{str} and Proposition \ref{ch} fail in the case of
a higher $\mathbb{Q}$-rank locally symmetric space.

\section{Positivity of the simplicial volume}

Now, we will prove the positivity of the simplicial volume of
$\mathbb{Q}$-rank one locally symmetric spaces $M$ covered by the product of
$\mathbb{R}$-rank one symmetric spaces.

Let $X_i$ be a complete, simply-connected, $n_i$-dimensional, Riemannian manifold with negative sectional
curvature bounded away from zero for each $i=1,\ldots,k$.
Let $X$ be the  $n$-dimensional product manifold $X_1 \times \cdots \times X_k$.
Now, we prove that the volume of geodesic $n$-simplices in $X$ is uniformly bounded
from above.

\begin{lemma}\label{vg}
The volume of geodesic $n$-simplices in $X$ is uniformly bounded from above.
\end{lemma}

\begin{proof}
Let $[x_0,\ldots,x_n]$ be  a geodesic $n$-simplex for an ordered
vertex set $\{ x_0,\ldots,x_n\} \subset X$. Let $p_i : X_1 \times \cdots \times
X_k \rightarrow X_i$ denote the projection map from $X$ onto $X_i$
for each $i=1,\ldots,k$. Then, we have
\begin{eqnarray}\label{voleqn}
\text{Vol}([x_0,\ldots,x_n]) \leq \prod_{i=1}^k \text{Vol}(p_i[x_0,\ldots,x_n]).
\end{eqnarray}

From Equation (\ref{voleqn}), it suffices to show that the volume of
$p_i[x_0,\ldots,x_n]$ in $X_i$ is uniformly bounded from above for
each $i=1,\ldots,k$. First, note that
$p_i([x_0,\ldots,x_n])=[p_i(x_0),\ldots,p_i(x_n)]$ because a
geodesic in $X$ projects to a geodesic in $X_i$ by the projection
map $p_i :X \rightarrow X_i$. In other words,
$p_i([x_0,\ldots,x_n])$ is a geodesic $n$-simplex in $X_i$. When $n
\geq n_i$,  a geodesic $n$-simplex in $X_i$ consists of at most
$\binom{n}{n_i}$ geodesic $n_i$-simplices in $X_i$. More precisely,
for a geodesic $n$-simplex $[y_0,\ldots,y_n]$ in $X_i$, we have
$$[y_0,\ldots,y_n] = \bigcup_{0\leq l_0 < \cdots <l_{n_i} \leq n}
[y_{l_0},\ldots,y_{l_{n_i}}].$$

Since the volume of geodesic $n_i$-simplices in $X_i$ is uniformly bounded from above \cite{IY82},
so is the volume of geodesic $n$-simplices in $X_i$.
Each uniform upper bound on the volumes of geodesic $n$-simplices in $X_i$
for each $i=1,\ldots,k$ gives a uniform upper bound on the volumes of geodesic $n$-simplices in $X$ by Equation (\ref{voleqn}),
which completes the proof.
\end{proof}

We now prove the main theorems. Let $M$ be an $n$-dimensional
Riemannian manifold. Then, the evaluation map
$$\langle \cdot , \cdot \rangle : C^*(M) \otimes C_*(M)\longrightarrow \mathbb{R}$$
is well-defined and it induces the Kronecker product on $H^*(M)\otimes
H_*(M)$. Let $K\subset M$ be a compact, connected subset with non-empty
interior. Let $\Omega^*(M,M-K)$ be the kernel of the restriction
homomorphism $\Omega^*(M)\rightarrow\Omega^*(M-K)$ on differential
forms. The corresponding cohomology groups are denoted by
$H^*_{dR}(M,M-K)$. The de Rham map $\Omega^*(M)\rightarrow C^*(M)$
restricts to the respective kernels and, thus, induces a homomorphism, called
relative de Rham map,
$$ \Psi^* :H^*_{dR}(M,M-K)\rightarrow H^*(M,M-K).$$

The relative de Rham map is an isomorphism. Note that integration gives a
homomorphism $\int :H^*_{dR}(M,M-K)\rightarrow\mathbb{R}$. Moreover, it
is well-known that
$$\langle\Psi^n[\omega],[M,M-K]\rangle=\int_M\omega$$ holds for all $n$-forms
$\omega$ \cite{Du76}.

\begin{proposition} \label{fc}
Let $M$ be a $n$-dimensional, $\mathbb{Q}$-rank one locally symmetric space covered by the
product of $\mathbb{R}$-rank one symmetric spaces. Let $c=\sum_{k\in\mathbb{N}}
a_k\sigma_k$ be a fundamental cycle of $M$ with $\| c \| _1 <\infty$. Then,
$$\sum_{k\in\mathbb{N}}a_k \cdot \langle dvol_M,st_n(\sigma_k)\rangle = vol(M).$$
\end{proposition}

\begin{proof}
Geodesic simplices in $M$ are
$C^1$, and the volume of geodesic simplices in $M$ is uniformly bounded from above
by Lemma \ref{vg}.
Hence, there exists a uniform constant $C>0$ such that we have
$$|\langle dvol_M,st_n(\sigma)\rangle|\leq C,$$
for any singular simplex $\sigma$ in $M$.
From this inequality, one can see that $\sum_{k\in\mathbb{N}}a_k \cdot \langle
dvol_M,st_n(\sigma_k)\rangle$ converges absolutely. Since $st^\text{lf}_*$ is chain
homotopic to the identity, $[st^\text{lf}_n(c)]$ is also a fundamental class of $M$.

Let $K\subset M$ be a connected, compact subset with non-empty interior.
For $\delta\in \mathbb{R}_{>0}$, let $g_\delta:M\rightarrow [0,1]$ be a
smooth function supported on the closed $\delta$-neighborhood $K_\delta$
of $K$ with $g_\delta |_K=1$. Then, $g_\delta dvol_M
\in\Omega^n(M,M-K_\delta)$ is a cocyle and
$$vol(K)=\lim_{\delta\rightarrow 0} \int_M g_\delta dvol_M.$$
The map $H_n(i_\delta):H^\text{lf}_n(M)\rightarrow H_n(M,M-K_\delta)$
induced by the inclusion $i_\delta:(M,\phi)\rightarrow(M,M-K_\delta)$ maps
the fundamental class of $M$ to the relative fundamental class
$[M,M-K_\delta ]$ of $(M,M-K_\delta)$, and $H_n(i_\delta)[st^\text{lf}_n(c)]$
is represented by $\sum_{im \sigma_k \cap K_\delta \neq \phi}a_k
st_n(\sigma_k)$. Since $H_n(i_\delta)[st^\text{lf}_n(c)]$ is also the relative
fundamental class of $(M,M-K_\delta)$, we have
\begin{eqnarray*}
\lefteqn{\lim_{\delta\rightarrow 0} \sum_{im \sigma_k \cap K_\delta \neq \phi}
a_k \cdot \langle g_\delta dvol_M, st_n(\sigma_k)\rangle} \\
&=&\lim_{\delta\rightarrow 0}\langle \Psi^n[g_\delta dvol_M], [M,M-K_\delta]\rangle \\
&=&\lim_{\delta\rightarrow 0} \int_M g_\delta dvol_M \\
&=& vol(K).
\end{eqnarray*}
For each $k\in\mathbb{N}$ and $\delta\in \mathbb{R}_{>0}$, we also have a
uniform upper bound $$|\langle g_\delta dvol_M, st_n(\sigma_k)\rangle
|\leq C,$$ and hence,
\begin{eqnarray*}
\Big| \sum_{k\in \mathbb{N}} a_k \cdot \langle dvol_M,st_n(\sigma_k) \rangle
-\sum_{im\sigma_k\cap K_\delta\neq\phi} a_k \cdot \langle g_\delta dvol_M,
st_n(\sigma_k) \rangle \Big| \\
\leq 2C \sum_{im\sigma_k \subset M-K} |a_k|.
\end{eqnarray*}

Because $\sum_{k\in\mathbb{N}}|a_k|<\infty$, there is an exhausting
sequence $(K^m)_{m\in \mathbb{N}}$ of compact, connected subsets of $M$
with non-empty interior satisfying
$$ \lim_{m\rightarrow\infty} vol(K^m)=vol(M)\text{ and }
\lim_{m\rightarrow\infty} \sum_{im\sigma_k \subset M-K^m} |a_k|=0.$$
Thus, the estimates of the previous paragraphs yield

{\setlength\arraycolsep{2pt}
\begin{eqnarray*}
\sum_{k\in \mathbb{N}} a_k \cdot \langle dvol_M,st_n(\sigma_k)\rangle &=&
\lim_{m\rightarrow\infty}\lim_{\delta\rightarrow\infty}
\sum_{im\sigma_k\cap K^m_\delta \neq \phi} a_k \cdot \langle g^m_\delta
dvol_M,st_n(\sigma_k)\rangle \\
&=&\lim_{m\rightarrow\infty} vol(K^m)\\
&=& vol(M),
\end{eqnarray*}}
which establishes the formula.
\end{proof}

From Proposition \ref{fc}, the positivity of the simplicial volume of
$\mathbb{Q}$-rank one locally symmetric spaces covered by
the product of $\mathbb{R}$-rank one symmetric spaces is directly obtained as follows.

\begin{theorem} \label{mt2}
Let $M$ be a $\mathbb{Q}$-rank one locally symmetric space covered by the product of
$\mathbb{R}$-rank one symmetric spaces.
Then, the simplicial volume of $M$ is positive.
\end{theorem}

\begin{proof}
Let $c=\sum_{k\in\mathbb{N}} a_k\sigma_k \in C^\text{lf}_n(M)$ be a
fundamental cycle. From Proposition \ref{fc}, we have
$$vol(M)=\sum_{k\in\mathbb{N}}a_k \cdot \langle dvol_M,st_n(\sigma_k)\rangle.$$
A uniform upper bound of the volume of geodesic simplices in $M$ yields the inequality
$$vol(M) \leq C \cdot \sum |a_k |,$$
where $C>0$ depends only on the universal cover of $M$. Dividing and
passing to the infimum over all fundamental cycles, it provides the positive
lower bound
$$ \| M \| \geq vol(M)/ C >0.$$
Therefore, we conclude that the simplicial volume of $M$ is positive.
\end{proof}

Gromov \cite{Gro82} provided a lower bound
   for the minimal volume minvol$(M)$, which is defined as the infimum of
    volumes over all complete Riemannian metrics on $M$
   with sectional curvatures bounded between $-1$ and $1$ in terms
   of the simplicial volume of an
     $n$-dimensional smooth manifolds $M$:
   \begin{eqnarray}\label{minvol}
  \| M\| \leq (n-1)^n n!\cdot \text{minvol}(M).
   \end{eqnarray}
By Inequality (\ref{minvol}) and Theorem \ref{mt2}, we have the following corollary.

\begin{corollary}\label{minvolume}
 The minimal volume of
$\mathbb{Q}$-rank one locally symmetric spaces covered by the product of
$\mathbb{R}$-rank one symmetric spaces is positive.
\end{corollary}
See \cite{CF1, CF2} for different approach using Lipschitz maps.
\section{Degree theorem}

For any proper map
$$f:N \rightarrow M$$ between finite volume Riemannian manifolds, a
locally finite fundamental cycle is mapped to a locally finite cycle. Hence, usual
inequality
\begin{eqnarray}\label{degree1}
deg(f)\cdot ||M||\leq ||N||
\end{eqnarray}
holds.
\begin{theorem}\label{degree}
Let $N$ be a Riemannian $n$-dimensional manifold of finite volume with Ricci curvature
bounded below by $-(n-1)$ and $M$ be a $\mathbb{Q}$-rank one locally
symmetric space covered by the product of $\mathbb{R}$-rank one symmetric spaces.
For any proper map $f:N\rightarrow
M$ we have
$$deg(f)\leq C_n \frac{vol(N)}{vol(M)},$$ where $C_n$ depends only on $n$.
\end{theorem}
\begin{proof}By Gromov \cite{Gro82},
$$||N||\leq (n-1)^n n! vol(N).$$
Since
$$vol(M)\leq C\cdot ||M||,$$ from Equation (\ref{degree1}), we
get
$$  \frac{deg(f)}{C}vol(M)\leq deg(f)||M||\leq ||N||\leq (n-1)^n n!
vol(N).$$  Since, for a given dimension $n$, there are only finitely many
symmetric spaces, the constant $C$ depends only on $n$.
\end{proof}
Such a degree theorem is known by \cite{CF2, LS09-2} for proper Lipschitz
map $f$ with the sectional curvature of $N$ bounded above by 1 and
any $n$-dimensional locally symmetric manifold $M$ of finite volume.
Note that they  obtained the degree theorem for proper
Lipschitz map $f$ by verifying the positivity of the Lipschitz
simplicial volume of $M$. Our result about the positivity of the
ordinary simplicial volume of $M$ yields the degree theorem without
any Lipschitz condition on map $f$.

% This is the end of the last section.

% Finally we create the bibliography or list of references.

% Every LaTeX document must end with \end{document}.
\end{document}